\documentclass[10pt]{paper}
\usepackage{amssymb,latexsym, amsmath}
\usepackage{amscd}

\usepackage{colortbl}
\usepackage[table]{xcolor}
\usepackage{mdwtab}

\newtheorem{exm}{Example}
\newtheorem{rem}{Remark}
\newtheorem{theorem}{Theorem}
\newtheorem{lemma}{Lemma}

\newcommand{\R}{\mathbb{R}}

\newcommand{\ddim}{\mathrm{ddim\;}}
\newcommand{\dind}{\mathrm{dind\;}}

\newcommand{\g}{\mathfrak{g}}

\newcommand{\ad}{\mathrm{ad}}
\DeclareMathOperator{\Ad}{\mathrm{Ad}}

\DeclareMathOperator{\rank}{\mathrm{rank}}
 \DeclareMathOperator{\pr}{\mathrm{pr}}
\DeclareMathOperator{\tr}{\mathrm{tr}}

\newlength{\arrayrulewidthOriginal}
    \newcommand{\Cline}[2]{%
      \noalign{\global\setlength{\arrayrulewidthOriginal}{\arrayrulewidth}}%
      \noalign{\global\setlength{\arrayrulewidth}{#1}}\cline{#2}%
      \noalign{\global\setlength{\arrayrulewidth}{\arrayrulewidthOriginal}}}
\newcommand{\dst}{\displaystyle}

\begin{document}

\title{Geodesic flows on Riemannian g.o. spaces\footnote{Journal ref: Regular and Chaotic Dynamics {\bf 16} (2011), No. 5, 504–-513 ---
issue: the conference {\it Geometry, Dynamics, Integrable Systems}
GDIS 2010 (Serbia, September 7–-13, 2010) dedicated to the 60th
birthdays of B. A. Dubrovin, V. V. Kozlov, I. M. Krichever, and A.
I. Neishtadt}}

\author{ Bo\v zidar Jovanovi\'c\\
\small \rm Mathematical Institute SANU, Serbian Academy of
Sciences and Arts\footnote{Kneza Mihaila 36, 11000 Belgrade,
Serbia, e-mail: bozaj@mi.sanu.ac.rs}}

\maketitle

\begin{abstract}
We prove the integrability of geodesic flows on the Riemannian
g.o. spaces of compact Lie groups, as well as on a related class
of Riemannian homogeneous spaces having an additional principal
bundle structure.\footnote{{\bf MSC:} 70H06, 37J35, 53D25}
\end{abstract}

\section{Introduction}

\paragraph{1.1.}
Let $(Q,ds^2)$ be a Riemannian manifold. A geodesic in $Q$ is
called {\it homogeneous} if it is the orbit of an one-parameter
group of isometries of $Q$. A Riemannian homogeneous space
$(Q=G/H,ds^2)$ is a {\it g.o. (geodesic orbit) space}  with
respect to $G$ if every geodesic is the orbit of an one-parameter
subgroup of $G$. Every naturally reductive space, in particular a
homogeneous space $G/H$ with a normal metric $ds^2_0$, is a g.o.
space. The converse is not true and the first counter example was
found by Kaplan \cite{Ka}, which motivated the study of g.o.
spaces, e.g., see \cite{KoVa, Go, BKV, KoNi, DKN, Ta, AlAr, AlNi}.

In Lagrangian and Hamiltonian mechanics homogeneous geodesics are
related to the notion of {\it relative equilibria}, e.g, see
\cite{La, To} and references therein. Following the line of
Hamiltonian mechanics, in this paper we study geodesic flows on
g.o. spaces $(Q=G/H, ds^2)$ of compact Lie groups $G$. Since the
geodesics are orbits, it is naturally to expect that the
corresponding geodesic flows, as in the case of normal metrics
\cite{BJ1, BJ2, Jo0}, are completely integrable. Indeed, we prove
their integrability in a noncommutative sense \cite{N, MF2, LMV}.

\paragraph{1.2.}
Let $(\mathcal F_G,\{\cdot,\cdot\})$ be the algebra of
$G$-invariant analytic, polynomial in momenta functions on
$T^*(G/H)$, where $T^*Q$ is endowed with the canonical Poisson
bracket. The Hamiltonian of a normal metric $ds^2_0$ is a central
function within the algebra $(\mathcal F_G,\{\cdot,\cdot\})$. The
main observation is that this property can be used for a
characterization of g.o. metrics on $G/H$: the Hamiltonian of a
$G$-invariant metric $ds^2$ is a central function in $\mathcal
F_G$ if and only if $ds^2$ is a g.o. metric (Lemma 2, Section 2).

Now, the integrability of geodesic flows on g.o. spaces follows
from a general geometric construction given in \cite{BJ1, BJ2}.
Let $\Phi: T^*Q\to \g^*$ be the momentum mapping of the natural
$G$-action. According to the Noether theorem, the momentum mapping
$\Phi$ is conserved along the flows of $G$-invariant Hamiltonians:
$\{\mathcal F_G,\Phi^*(\R[\mathfrak g])\}=0$. Moreover, the
collection of functions
\begin{equation}\label{int}
\mathcal F=\mathcal F_G+\Phi^*(\R[\mathfrak g^*])
\end{equation}
is a {\it complete set} of functions on $T^*Q$:
\begin{equation}\label{potpun}
\ddim{\mathcal F}+\dind{\mathcal F}=\dim T^*Q.
\end{equation}
Here, for a set of functions $\mathcal F$ closed under the Poisson
bracket, such that the dimensions of $\mathrm F_x=\langle df(x)
\,\vert\, f\in\mathcal F\rangle$ and $\ker \{\cdot,\cdot\}
\vert_{\mathrm F_x\times\mathrm F_x}$ are constant on an open
dense set $U\subset M$, we denote $\ddim\mathcal F=\dim\mathrm
F_x$ and $\dind\mathcal F=\dim\ker \{\cdot,\cdot\} \vert_{\mathrm
F_x\times\mathrm F_x}$, $x\in U$ ({\it differential dimension} and
{\it differential index} of $\mathcal F$) \cite{BJ2, Jo0}.

Since central functions in $\mathcal F_G$ Poisson commute both
with $\mathcal F_G$ and $\Phi^*(\R[\mathfrak g^*])$,  we get:

\begin{theorem}
The geodesic flow of a g.o. space $(Q=G/H, ds^2)$ is completely
integrable with a complete noncommutative set of integrals
\eqref{int}.
\end{theorem}

\paragraph{1.3.}
Mishchenko and Fomenko stated the conjecture that noncommutative
integrability implies the usual Liouville one by means of
integrals that belong to the same functional class as
noncommutative integrals \cite{MF2, Sa, BJ2}. One can always
construct a complete commutative subset within $\Phi^*(\mathfrak
g^*)$. Thus, for a complete set \eqref{int}, the conjecture
reduces to a construction of a complete commutative subset
$\mathcal B \subset \mathcal F_G$ (see \cite{BJ2, BJ3})
\begin{equation}
\ddim {\mathcal B} = \dind\mathcal B=\frac12\left(\ddim \mathcal
F_G + \dind \mathcal F_G\right). \label{C}
\end{equation}

If the required subset $\mathcal B\subset \mathcal F_G$ exist, we
say that $(G,H)$ is an {\it integrable pair}. There are several
known classes of integrable pairs (see \cite{Th, Mik, MS, BJ1,
BJ3, MP, DGJ, Jo}) but the general problem rest still unsolved.

We recall the known examples of g.o. spaces $G/H$ with non-normal
metrics of compact Lie groups (Section 2). In most of them, the
underlying metrics are deformations $ds^2_\lambda$ of a normal
metric in the fiber direction of the  bundle
\begin{equation}\label{bundle}
K/H \longrightarrow Q=G/H \longrightarrow M=G/K,
\end{equation}
where $K$ is a subgroup of $G$, $H\subset K\subset G$.

Tamaru classified triples $(G,K,H)$ where $(G,K)$ is a compact
effective irreducible symmetric pair and $(G/H,ds^2_\lambda)$ is a
g.o. space \cite{Ta}. It turns out that these homogeneous spaces
have the complexity equal 0 or 1 (see the classification of
complexity 0 and 1 homogeneous spaces of the reductive Lie groups
\cite{Mik, Ya, Pa, MS, ArCh}). This implies, by the result of
Mykytyuk and Stepin \cite{MS}, that the appropriate geodesic flows
are Liouville integrable by means of analytic functions,
polynomial in momenta (Theorem 4).
%The classification of complexity 0 and 1 homogeneous spaces of the
%reductive Lie groups is well known (see \cite{Mik, Ya, Pa, MS,
%ArCh}), thus the considered g.o. spaces with non-normal metrics do
%not provide new integrable pairs.
For a reader's sake, by
analyzing the algebra of $G$-invariant functions $\mathcal F_G$ a
direct derivation of the complexity and commuting integrals for
various g.o. spaces is given (Theorem 4, Examples 4 -- 7, Section
4).

As a bi-product, in Section 3 we obtain a class of Riemannian
homogeneous spaces, defined by fibration \eqref{bundle}, with
completely integrable geodesic flows (Theorem \ref{pomoc},
Examples 1 -- 3). In particular, we have (see Example 1):

\begin{theorem}
The homogeneous spaces fibered over irreducible symmetric spaces
with invariant Einstein metrics classified by Jensen \cite{Je}
have completely integrable geodesic flows.
\end{theorem}

\section{Homogeneous g.o. spaces}

\paragraph{2.1. Homogeneous Riemannian spaces.}
A Riemannian manifold $(Q,ds^2)$ is homogeneous if it admits a
transitive connected Lie group of isometries. Then $Q$ can be seen
as a coset space $G/H$ with a $G$-invariant metric, where $H$ is
the isotropy group at some point $o\in Q$. Since $H$ is compact,
there exist an $\Ad_H$-invariant subspace $\mathfrak v\subset
\mathfrak g$ ($[\mathfrak h,\mathfrak v]\subset \mathfrak v$) such
that $\mathfrak g=\mathfrak h\oplus \mathfrak v$ (a reductive
decomposition of $\mathfrak g$). Here $\mathfrak g$ and $\mathfrak
h$ are Lie algebras of $G$ and $H$ and $\Ad_H$ is the restriction
of the adjoint representation $\Ad_G$ to $H$.

We can identify $\mathfrak v$ with $T_o Q$ by taking the value of
the corresponding Killing vector field at $o$. In this way, the
isotropy representation of $H$ at $T_o Q$ is identified with the
restriction of the adjoint representation of $H$ to $\mathfrak v$.
The restriction of a $G$-invariant metrics $ds^2$ to $T_o Q \cong
\mathfrak v$ defines $\Ad_H$-invariant scalar product
$(\cdot,\cdot)$ on $\mathfrak v$. Conversely, for each
$\Ad_H$-invariant scalar product $(\cdot,\cdot)$ on $\mathfrak v$,
there exist an unique $G$-invariant Riemannian metric $ds^2$ on
$Q=G/H$ such that its restriction to $T_o Q\cong \mathfrak v$ is
given by $(\cdot,\cdot)$ (e.g., see \cite{Be}).

A $G$-invariant metric $ds^2_0$ is called {\it normal} if
$(\cdot,\cdot)$ is the restriction of an $\Ad_G$-invariant
non-degenerate symmetric bilinear form $q$ on $\mathfrak g$. If
$G$ is semi-simple and $q=-B$ ($B$ is the Killing form), the
metric $ds^2_0$ is called {\it standard}.

A metric $ds^2$ is {\it naturally reductive} if
\begin{equation}\label{nr}
([X,Y]_\mathfrak v,Z) +(Y,[X,Z]_\mathfrak v)=0, \qquad
X,Y,Z\in\mathfrak v.
\end{equation}

Clearly, a normal metric is naturally reductive. A notion of
natural reductivity depends on the choice of the group $G$ and the
choice of $\mathfrak v$. A result of Konstant (e.g., see
\cite{Be}, page 196) stated that $(Q=G/H,ds^2)$ is naturally
reductive if and only if there exist a subgroup $\tilde G \subset
G$ with the Lie algebra $\tilde{\mathfrak g}=\mathfrak
v+[\mathfrak v,\mathfrak v]$, such that $ds^2$ is a normal metric
on $Q=\tilde G/(\tilde G \cap H)$.

\paragraph{2.2. Riemannian g.o. spaces}
The condition \eqref{nr} is equivalent to the following
geometrical property: the curve $\gamma(t)=\exp(tX)\cdot o$ is a
geodesic for all $X\in\mathfrak v$. Thus, for a naturally
reductive homogeneous space every geodesic is homogeneous. By
generalizing the property above, Kowalski and Vanhecke \cite{KoVa}
introduced a notion of g.o. spaces. A Riemannian homogeneous space
$(Q=G/H,ds^2)$ is a {\it g.o. space}, if every geodesic is an
orbit of an one-parameter subgroup $\exp(t X)$, $X\in\mathfrak g$.

The following algebraic  condition on an orbit
 to be a geodesic is well known \cite{KoVa}:

\begin{lemma}[Geodesic lemma]
Let $X\in\mathfrak v$, $F\in\mathfrak h$. The orbit
$\gamma(t)=\exp(t(X+F))\cdot o$ through the point $o=[H]\in G/H$
is a geodesic of $(Q=G/H,ds^2)$ if and only if
\begin{equation}\label{gl}
([X+F,Y]_\mathfrak v,X)=0, \qquad\text{for all} \qquad
Y\in\mathfrak v.
\end{equation}
\end{lemma}

Therefore, a Riemannian homogeneous space $(Q=G/H,ds^2)$ is a g.o.
space if and only if  for any vector $X\in\mathfrak v$ there exist
$F=F(X)\in\mathfrak h$ that satisfies \eqref{gl}.

Examples of such spaces are symmetric spaces, weakly symmetric
spaces (for any two points $p,q \in Q$, there exist an isometry
which interchanges $p$ and $q$ \cite{Se, BeVa, BKV,  Vi, Ya, Zi2})
and naturally reductive spaces.

\paragraph{2.3.  $G$-invariant functions and a characterization of g.o. metrics.} From now on $G$ is a compact
connected Lie group $G$. We fix an $\Ad_G$-invariant scalar
product $\langle \cdot,\cdot \rangle$ on $\mathfrak g$. Then
$\mathfrak g=\mathfrak h\oplus\mathfrak v$, $\mathfrak v=\mathfrak
h^\perp$ is a reductive decomposition of $\mathfrak g$.  An
$\Ad_H$-invariant scalar  product on $\mathfrak v\cong T_o Q$ can
be represented as
$$
(\,\cdot\,,\,\cdot\,)=\langle I(\,\cdot\,),\,\cdot\, \rangle,
$$
where $I: \mathfrak v\to \mathfrak v$ is a positive definite,
$\Ad_H$-invariant operator. For a given operator $I$, we denote
the associated $G$-invariant metric by $ds^2_I$. In particular,
for $I$ equals to the identity, we have a normal metric $ds^2_0$.

By using $ds^2_0$ we identify tangent and cotangent bundles $TQ
\cong T^*Q$. Within this identification, the $G$-invariant
functions are in one-to-one correspondence with $\Ad_H$-invariant
polynomials $\R[\mathfrak v]^H$ on $\mathfrak v$.  The Hamiltonian
of the metric $ds^2_{I}$ is given by $h_A(x)=\frac12\langle A
x,x\rangle$ where $A=I^{-1}$, while the Hamiltonian of the normal
metric $ds^2_0$ is simply $h_0(x)=\frac12\langle x,x\rangle$.

Further, the canonical Poisson bracket on $T^* Q$ corresponds to
the restriction of the Lie-Poisson bracket to
$\mathbb{R}[\mathfrak v]^{H}$ (see \cite{Th}):
\begin{equation}
\{f,g\}_\mathfrak v(x)=-\langle x,[\nabla f(x),\nabla
g(x)]\rangle, \qquad f,g: \mathfrak v\to \R.\label{rpz}
\end{equation}

Let $\g_x$ and $\mathfrak h_x$ be the isotropy algebras of $x$  in
$\g$ and $\mathfrak h$ and let $\mathfrak j_x \subset \mathfrak v$
be the orthogonal complement to the orbit
$\Ad_H(x)\subset\mathfrak v$: $\mathfrak j_x=\{ \eta\in \mathfrak
v \, \vert\, \langle \eta,[x,\mathfrak h]\rangle=0\}$.  For a
generic point $x\in \mathfrak v$, it is spanned by gradients of
polynomials in $\mathbb{R}[\mathfrak v]^{H}$. The maximal number
of functionally independent polynomials is equal to
\begin{equation}\label{ddim}
\ddim\mathcal F_G=\ddim\R[\mathfrak v]^H=\dim\mathfrak j_x=\dim
\mathfrak v-\dim \mathfrak h +\dim \mathfrak h_x,
\end{equation}
for a generic $x\in\mathfrak v$.

An $\Ad_H$-invariant polynomial $f(x)$ is {\it central (Casimir)}
within $\R[\mathfrak v]^H$ if it commutes with all invariant
polynomials, that is $\nabla f(x)\in \ker\Lambda_x$, where
\begin{equation*} \label{LPB}
\Lambda_x(\eta_1,\eta_2)=-\langle x,[\eta_1,\eta_2]\rangle, \qquad
\eta_1,\eta_2 \in \mathfrak j_x.
\end{equation*}

The kernel of $\Lambda_x$ is $ \ker\Lambda_x=\pr_\mathfrak
v\g_x\subset \mathfrak j_x$. Thus, the maximal number of
functionally independent central polynomials is $\dind\R[\mathfrak
v]^H=\dim\ker\Lambda_x=\dim \mathfrak g_x-\dim\mathfrak h_x,$ for
a generic $x\in\mathfrak v$ \cite{BJ2, BJ3}. Moreover, since
$[\g_x,\g_x]\subset\mathfrak h_x$ (see Mykytyuk \cite{Mik}) we
have also:
\begin{equation}\label{dind}
 \dind\mathcal
F_G=\dind\R[\mathfrak v]^H=\dim \mathfrak g_x-\dim\mathfrak
h_x=\rank\g-\rank\mathfrak h_x,
\end{equation}
for a generic $x\in\mathfrak v$.

Let $f_1,\dots,f_r$ be a base of homogeneous $\Ad_G$-invariant
polynomials on $\g$ ($r=\rank\g$). They are central functions with
respect to the Lie-Poisson brackets on $\g$, while their
restrictions $p_i=f_i\vert_\mathfrak v$ to $\mathfrak v$ are
central functions in $\R[\mathfrak v]^H$. Among $p_1,\dots,p_r$
there are $\dind\R[\mathfrak v]^H$ functionally independent
polynomials. In particular, a normal metric Hamiltonian
$h_0(x)=\frac12\langle x,x\rangle$ is a central function.

\begin{lemma}\label{central}
A Riemannian homogeneous space $(Q=G/H, ds^2_I)$ is a g.o. space
if and only if the Hamiltonian $h_A(x)=\frac12\langle A
x,x\rangle$ is a central function in $\R[\mathfrak v]^H$.
\end{lemma}

\noindent{\it Proof.} The Hamiltonian $h_A(x)=\frac12\langle A
x,x\rangle$ is $\Ad_H$-invariant implying
\begin{equation}\label{cond}
\ad_\mathfrak h(x)\cdot h_A(x) = 0 \quad \Longleftrightarrow\quad
[x,Ax]_\mathfrak h=0.
\end{equation}

According to Lemma 1, $(Q=G/H, ds^2_I)$ is a g.o. space if and
only if for every $X\in\mathfrak v$ there exist $F(X)\in\mathfrak
h$ satisfying
\begin{equation}\label{go}
\langle [X+F(X),Y]_\mathfrak v,IX \rangle=0, \qquad\text{for all}
\qquad Y\in\mathfrak v.
\end{equation}

Denote $x=IX$, $a=a(x)=F(Ax)$. Then we can rewrite \eqref{go} as
\begin{equation}\label{go2}
[a(x)+Ax,x]_\mathfrak v=0.
\end{equation}

Since $[a,x]\in\mathfrak v$, the relation \eqref{cond} implies
equivalence of \eqref{go2} and
\begin{equation*}\label{go3}
 [a(x)+Ax,x]=0 \quad \Longleftrightarrow \quad Ax=\nabla
h_A(x)\in \pr_\mathfrak v\mathfrak g_x=\ker\Lambda_x.
\end{equation*}

Therefore $(Q=G/H, ds^2_I)$ is a g.o. space if and only if
 $h_A(x)$ is a central function in $\R[\mathfrak v]^H$. \hfill$\Box$

\begin{rem}{\rm
The geodesic flows of a g.o. metric $ds^2_I$ and of a normal
metric $ds^2_0$ share the same invariant isotropic toric foliation
 defined by integrals \eqref{int} (a generic torus has a dimension $\dind\mathcal F=\dind\mathcal F_G=\dind
\Phi^*(\R[\mathfrak g^*])$). However, the tori can be resonant and
the closures of geodesic lines might be different. For example,
while the geodesic lines of the standard metric on a sphere are
great circles, we see that generic geodesic lines on the distance
spheres $S^{4n+3}=Sp(n+1)/Sp(n)$ and $S^{2n+1}=SU(n+1)/SU(n)$,
described in Examples \ref{sfera1} and \ref{sfera2}, filled up
regions that are projections of two-dimensional invariant
isotropic tori.}\end{rem}

\paragraph{2.4. Examples.}
Here we list some known examples of g.o. spaces with compact
groups of isometries and with non-normal metrics.

Aleksieevsky and Arvanitoyeorgos proved that among all flag
manifolds $Q=G/H$ of simple Lie groups $G$ only the manifolds
\begin{equation}\label{flag}
SO(2n+1)/U(n), \qquad Sp(n)/U(1)\times Sp(n-1)
\end{equation}
admit invariant metrics with homogeneous geodesics, not homothetic
to the standard metric. These manifolds have an one-parameter
family of metrics $ds^2_\lambda$, $\lambda>0$ with homogeneous
geodesics and are all weakly symmetric, see \cite{Zi2}. The metric
$ds^2_1$ is the standard one. It has the full connected isometry
group $SO(2n+2)$ (respectively $SU(2n-1)$) and is the standard
metric of the symmetric space $SO(2n+2)/U(n+1)$ (respectively the
complex projective space $\mathbb{CP}^{2n-1}=SU(2n-1)/U(2n-2)$).
All the other metrics have the full connected isometry group
$SO(2n+1)$ (respectively $Sp(n)$) and the corresponding spaces are
not naturally reductive. This two classes of examples exhaust all
simply connected compact irreducible Riemannian non-normal g.o.
manifolds of positive Euler characteristic (see Aleksieevsky and
Nikonorov \cite{AlNi}).

For $l=2$, $Sp(2)/U(1)\cdot Sp(1) \cong SO(5)/ U(2)$ is the
6-dimensional not naturally reductive g.o. space of Kowalski and
Vanhecke, who classified all g.o. spaces in dimension $\le 6$
\cite{KoVa}.

Du\v sek, Kowalski and Nik\v cevi\' c constructed a family of
invariant metrics $ds^2_{p,q}$ on a 7-dimensional homogeneous
space $Q=(SO(5)\times SO(2))/U(2)$, such that $(Q,ds^2_{p,q})$ are
g.o. spaces which are not naturally reductive \cite{DKN}. Besides,
the group $SO(5)$ acts as a transitive group of isometries, but
the homogeneous spaces $(Q=SO(5)/SU(2),ds^2_{p,q})$ are not g.o.
spaces. Note that $(SO(5)\times SO(2))/U(2)$ is also a weakly
symmetric space (see the classification of weakly symmetric spaces
given by Yakimova \cite{Ya}).

Finally, Tamaru classified  the Riemannian, non-standard g.o.
spaces, which are fibered over the irreducible symmetric spaces
\cite{Ta0, Ta} (see Table 1 below). In particular, the flag
manifolds \eqref{flag} belong to the considered class of
homogeneous spaces. They are associated to the triples
$(SO(2n+1),SO(2n),U(n))$ and $(Sp(n),Sp(1)\times
Sp(n-1),U(1)\times Sp(n-1))$, respectively. Note that the spaces
listed in Table 1 can be naturally reductive with respect to a
suitably larger symmetry group (e.g., see Kowalski and Nik\v
cevi\' c \cite{KoNi} for the space $SU(3)/SU(2)\cong S^5$).

\section{Geodesic flows on fiber bundles}

\paragraph{3.1. Integrable pairs and complexity of homogeneous spaces.}
Recall that for a complete set of integrals \eqref{int}, the
Mishchenko--Fomenko conjecture reduces to a construction of a
commutative subset $\mathcal B \subset \mathcal F_G$ that
satisfies condition \eqref{C}. If such a set exists, we call
$(G,H)$ an {\it integrable pair}. If $G/H$ is a weakly symmetric
space, in particular if $G/H$ is a symmetric space, the algebra
$\mathcal F_G$ is already commutative \cite{Vi}. In this case the
pair $(G,H)$ is obviously integrable and we need only Noether's
integrals to integrate the geodesic flow (see \cite{Br, Mik}).
Further examples can be found in \cite{Th, MS, BJ1, BJ3, MP, DGJ,
Jo}.

There is a well known notion of  homogeneous spaces {\it
complexity} of complex reductive Lie groups \cite{Pa}. In the case
of a homogeneous space $G/H$ of a compact Lie group $G$, the
complexity of $G^\mathbb{C}/H^\mathbb{C}$ corresponds to the
number of independent polynomials, apart from the central
polynomials, which we need to form a complete commutative subset
$\mathcal B \subset \R[\mathfrak v]^H$ (see Mykytyuk \cite{Mik2}).
We refer to
$$
c(G,H)=\frac12(\ddim\R[\mathfrak v]^H-\dind\R[\mathfrak v]^H),
$$
as a {\it complexity} of the homogeneous space $G/H$ or the pair
$(G,H)$.

Symmetric and weakly symmetric spaces have the complexity $c=0$
\cite{Vi}. The pairs $(G,H)$ with the complexity $c=1$ are also
integrable: we can take an arbitrary non central polynomial and
the central polynomials to form a required commutative set
$\mathcal B$ (see Mykytyuk and Stepin \cite{MS}). The
classification of the complexity 1 homogeneous spaces of reductive
algebraic groups is given by Arzhantsev and Chuvashova \cite{ArCh}
(see also \cite{Pa, MS}).

\paragraph{3.2. Triples of Lie groups.}
Suppose $G$ is a connected compact semi-simple Lie group. Let $H
\subset K \subset G$ be a chain of subgroups and $\mathfrak h
\subset \mathfrak k \subset \mathfrak g$ be the appropriate chain
of subalgebras.

Let $\langle \cdot,\cdot\rangle$ to be the negative of the Killing
form on $\mathfrak g$. Taking the suitable orthogonal complements,
we obtain decompositions
\begin{equation}\label{dec}
\mathfrak g=\mathfrak h\oplus\mathfrak v=\mathfrak h\oplus
\mathfrak l \oplus \mathfrak m, \qquad \mathfrak k=\mathfrak
h\oplus \mathfrak l.
\end{equation}

Following \cite{Ta0, Ta}, consider a deformation $ds^2_\lambda$ of
the standard metric on $Q=G/H$ in the direction of the fiber $K/H$
of the bundle \eqref{bundle}, defined as a scalar product on
$\mathfrak v=\mathfrak l \oplus \mathfrak m \cong T_oQ$:
\begin{equation}\label{lambda}
(\cdot,\cdot)_\lambda=\lambda \langle
\cdot,\cdot\rangle_{\mathfrak l \times \mathfrak l}+ \langle
\cdot,\cdot\rangle_{\mathfrak m \times \mathfrak m}, \qquad
\lambda>0.
\end{equation}

The Hamiltonian of the metric $ds^2_\lambda$ is a deformation of
the Hamiltonian of the standard metric $ds^2_0$:
$$
h_\lambda(x)=h_0(x)+(\lambda^{-1}-1)\delta(x), \quad
h_0(x)=\frac12\langle x_\mathfrak m,x_\mathfrak
m\rangle+\frac12\langle x_\mathfrak l,x_\mathfrak l\rangle, \quad
\delta(x)=\frac12\langle x_\mathfrak l,x_\mathfrak l\rangle,
$$
where $x=x_\mathfrak l+x_\mathfrak m$, $x_\mathfrak l\in \mathfrak
l$, $x_\mathfrak m\in\mathfrak m$.

A natural problem is whether the geodesic flow of  $ds^2_\lambda$
is integrable. The positive answer for a class of spaces is given
in \cite{BJ3}. Here we give the following simple criterium.

\begin{theorem}\label{pomoc}
Let $(G,K,H)$ be a triple of connected Lie groups such that
 $[\mathfrak l,\mathfrak h]=0$ and $[\mathfrak l,\mathfrak
 l]\subset \mathfrak l$. Then

(i) the geodesic flow of the metric $ds^2_\lambda$ is
noncommutative integrable;

(ii) if $\mathfrak l$ is commutative  ($\mathfrak l$ is a subset
of the center of $\mathfrak k$) and $(G,K)$ is an integrable pair
then $(G,H)$ is also an integrable pair.
\end{theorem}

\noindent{\it Proof.} (i) Let $L$ be the connected subgroup of $K$
with the Lie algebra $\mathfrak l$. Then $H$ and $L$ are normal
subgroups of $K$. Locally, $K$ is isomorphic to the product
$H\times L$. The right translations on $(Q=G/H,ds^2_\lambda)$ with
elements of $L$ are well defined isometries. Note that the actions
of $G$ and $L$ commute and that \eqref{bundle} is a $L$-principal
bundle.

Let $\Psi: T^*Q\to \mathfrak l^*$ be the momentum mapping of the
$L$-action lifted to the cotangent bundle and let $\mathcal
F_{G\times L}$ be the algebra of $G\times L$-invariant analytic,
polynomial in momenta functions on $T^*Q$. By the general
construction \cite{BJ2}, the collection of functions
\begin{equation}
\mathcal F_{G\times L}+\Psi^*(\R[\mathfrak
l^*])+\Phi^*(\R[\mathfrak g^*]) \label{int2}
\end{equation}
is a complete set on $T^*Q$. It remains to observe that these
functions commute both with $h_0(x)$ and $\delta(x)$.

\

(ii) Within the identification $\mathcal F_G\cong \R[\mathfrak
v]^H$, the algebra $\Psi^*(\R[\mathfrak l^*])$ corresponds to the
algebra $\mathcal L$ generated by linear functions on $\mathfrak
l$, extended to $\mathfrak v$ via decomposition \eqref{dec}.

Assume $\mathcal B_M$ is a complete commutative set in
$\R[\mathfrak m]^K$. Let $\mathcal B=\mathcal B_M+\mathcal L$,
where the polynomials in $\mathcal B_M$ are extended to $\mathfrak
v$ via decomposition \eqref{dec}. It is clear that $\mathcal B$ is
a commutative subset of $\R[\mathfrak v]^H$. Indeed, from
$[\mathfrak l,\mathfrak k]=0$ we get $\{\mathcal L,\mathcal
L\}_\mathfrak v=0$ and $\{\mathcal B_M,\mathcal L\}_\mathfrak v=0$
follows from $\Ad_K$-invariancy of polynomials in $\mathcal B_M$.
Furthermore, since
\begin{equation}\label{isotropy}
\dim \g_x \le \dim g_{x_\mathfrak m}, \quad \text{for a generic}
\quad x\in\mathfrak v, \, x_\mathfrak m\in\mathfrak m,
\end{equation}
we get
\begin{eqnarray*}
\ddim\mathcal B &=& \ddim \mathcal B_M+\dim\mathfrak l=
\frac12(\ddim\R[\mathfrak m]^K+\dind\R[\mathfrak
m]^K)+\dim\mathfrak l\\
&=& \frac12(\dim\mathfrak m-\dim \mathfrak k +\dim\mathfrak
g_{x_\mathfrak m})+\dim\mathfrak l= \frac12(\dim\mathfrak v-\dim\mathfrak h+\dim\g_{x_\mathfrak m})\\
&\ge &\frac12(\ddim\R[\mathfrak v]^H+\dind\R[\mathfrak v]^H).
\end{eqnarray*}

On the other hand, for a commutative set $\mathcal B$ we always
have the inequality $ \ddim\mathcal B \le
\frac12(\ddim\R[\mathfrak v]^H+\dind\R[\mathfrak v]^H). $
Therefore, \eqref{isotropy} is an equality and $\mathcal B$ is a
complete  set.  \hfill$\Box$

\begin{exm}{\rm
The homogeneous spaces considered in Theorem \ref{pomoc} are
studied from a point of view of a construction of invariant
Einstein metrics  \cite{Je} (see also \cite{Be}). Jensen
classified Einstein metrics of the form $ds^2_\lambda$ on
homogeneous spaces $Q=G/H$ fibered over irreducible symmetric
spaces $M=G/K$. When $\mathfrak l$ is nonabelian, then there are
two Einstein metrics, neither of them equals to the standard one.
When $\mathfrak l$ is commutative, there is only one Einstein
metric, which equals to the standard one when $\dim L=1$ and $\dim
M=2$ \cite{Je}. For example, the triple $(SO(n),SO(r)\times
SO(n-r),SO(n-r))$ provides two Einstein metrics on a Stiefel
variety $V_{n,r}=SO(n)/SO(n-r)$. The geometry and integrability of
natural mechanical systems on Stiefel varieties $V_{n,r}$ endowed
with metrics  $ds^2_\lambda$ and with quadratic potentials is
studied in \cite{FJ}.}\end{exm}

\begin{exm}{\rm A metric $ds^2_\lambda$ can be always additionally perturb by taking
an arbitrary integrable Euler equation on $\mathfrak l\cong
\mathfrak l^*$
$$
\dot x_\mathfrak l=[x_\mathfrak l,A x_\mathfrak l], \qquad
A:\mathfrak l \to \mathfrak l
$$
with a complete set of integrals $\mathcal S$. Namely, define $
h_{\lambda,A}=h_\lambda(x)+\frac12\langle x_\mathfrak l,A
x_\mathfrak l\rangle. $ We can find $\lambda$ such that
$h_{\lambda,A}$ is positive definite. It is a Hamiltonian of a
suitable $G$-invariant metric $ds^2_{\lambda,A}$. From the
completeness of sets $\mathcal S$ on $\mathfrak l$ and
\eqref{int2} on $T^*Q$, we get that $\mathcal F_{G\times
L}+\mathcal S$ is a complete subset of $\mathcal F_G\cong
\R[\mathfrak v]^H$:
\begin{equation}\label{int3}
\ddim{(\mathcal F_{G\times L}+\mathcal S)}+\dind{(\mathcal
F_{G\times L}+\mathcal S)}=\ddim\R[\mathfrak
v]^H+\dind\R[\mathfrak v]^H.
\end{equation}
Thus, the geodesic flow of $ds^2_{\lambda,A}$ is completely
integrable. Note that Einstein metrics on Stiefel manifolds
$V_{n,r}$ obtained in \cite{ADN} are of the form
$ds^2_{\lambda,A}$ \cite{FJ}. }\end{exm}

\begin{exm}{\rm As an application of Theorem \ref{pomoc}, item (ii), let us consider integrable
pairs $ (Sp(n),U(k_1)\times\dots\times U(k_r)) $ ($k_i\le
[(n+1)/2]$, see \cite{BJ3}). Then we get that
$(Sp(n),SU(k_1)\times\dots\times SU(k_r))$ are integrable pairs as
well.}\end{exm}

\section{G.O. spaces fibered over symmetric spaces}

\paragraph{4.1.}
Consider a triple of Lie groups $(G,K,H)$. According to Lemma
\ref{central}, $(G/H,ds^2_\lambda)$ is a g.o. space if and only if
$\delta(x)$ is a central function in $\R[\mathfrak v]^H$, that is
$\nabla\delta(x)=x_\mathfrak l\in\pr_{\mathfrak v}\g_x$. In such a
way we obtain the following condition (see Gordon \cite{Go}):
$(G/H,ds^2_\lambda)$ is a g.o. space if and only if for every
$x\in\mathfrak v$, there exist $a=a(x)\in\mathfrak h$, that
satisfies
$$
[a+x_\mathfrak l,x_\mathfrak l+x_\mathfrak k]=0 \quad
\Longleftrightarrow \quad [a,x_\mathfrak l]=0, \quad
[a+x_\mathfrak l,x_\mathfrak m]=0.
$$
In particular, if $(G/H,ds^2_\lambda)$ is a g.o. space for a
single $\lambda\ne 1$, then it is a g.o. space for all
$\lambda>0$. Triples $(G,K,H)$ where $(G,K)$ is a compact
effective irreducible symmetric pair and $(G/H,ds^2_\lambda)$ is a
g.o. space are classified by Tamaru \cite{Ta}.

\begin{table}[ht]
\begin{tabular}[C]{!{\vline[1.2pt]}c!{\vline[1.2pt]}c|c|c!{\vline[1.2pt]}c !{\vline[1.2pt]}}
\Cline{1.2pt}{1-5}  &  $\,\dst \quad\mathfrak g\;\quad$ & $\dst
\,\quad\mathfrak k\quad\;$ & $\dst \,\quad\mathfrak h\quad\;$ &     \\
\Cline{1.2pt}{1-5} 1 & $\dst \quad so(2n+1) \quad$ & $\dst \quad
so(2n) \quad$ & $\dst \quad u(n) \quad$ &
  $\dst \, n\ge 2 \,$  \\
\hline  2 & $\dst \quad so(4n+1) \quad$ & $\dst \quad so(4n)
\quad$ & $\dst \quad su(2n) \quad$ & $\dst \quad n\ge 1 \quad$  \\
\hline  3 & $\dst \quad so(8) \quad$ & $\dst \quad so(7) \quad$ &
$\dst \, g_2 \,$ &  \\
\hline  4 & $\dst \quad so(9) \quad$ & $\dst \quad so(8) \quad$ &
$\dst \quad spin(7) \quad$ &    \\
\hline  5 & $\dst \quad su(n+1) \quad$ & $\dst \quad u(n) \quad$ &
$\dst \quad su(n) \quad$ &
  $\dst \, n\ge 2 \,$  \\
\hline  6 & $\dst \quad su(2n+1) \quad$ & $\dst \quad u(2n) \quad$
& $\dst \quad u(1)\oplus sp(n) \quad$ &
  $\dst \, n\ge 2 \,$  \\
\hline  7 & $\dst \quad su(2n+1) \quad$ & $\dst \quad u(2n) \quad$
& $\dst \quad  sp(n) \quad$ &
  $\dst \, n\ge 2 \,$  \\
\hline  8 & $\dst \quad sp(n+1) \quad$ & $\dst \quad sp(1)\oplus
sp(n) \quad$ & $\dst \quad  u(1)\oplus sp(n) \quad$ &
  $\dst \, n\ge 1 \,$  \\
\hline  9 & $\dst \quad sp(n+1) \quad$ & $\dst \quad sp(1)\oplus
sp(n) \quad$ & $\dst \quad  sp(n) \quad$ &
  $\dst \, n\ge 1 \,$  \\
%\Cline{1.2pt}{1-5}
\hline 10 & $\dst \quad su(2r+n) \quad$ & $\dst \, su(r)\oplus
su(r+n)\oplus\R \,$ & $\dst \quad su(r)\oplus
su(r+n) \quad$ &  $\dst \, r\ge 2, n\ge 1 \,$  \\
\hline  11 & $\dst \quad so(4r+2) \quad$ & $\dst \quad u(2r+1)
\quad$ & $\dst \quad su(2r+1) \quad$ & $\dst \quad r\ge 2 \quad$  \\
\hline 12 & $\dst \quad e_6 \quad$ & $\dst \quad so(10)\oplus\R
\quad$ & $\dst \quad so(10)\quad$ &    \\
%\Cline{1.2pt}{1-5}
\hline13 & $\dst \quad so(9) \quad$ & $\dst \quad so(7)\oplus
so(2) \quad$ & $\dst \quad g_2\oplus
so(2) \quad$ &   \\
\hline  14 & $\dst \quad so(10) \quad$ & $\dst \quad so(8)\oplus
so(2)
\quad$ & $\dst \quad spin(7)\oplus so(2) \quad$ &   \\
\hline 15 & $\dst \quad so(11) \quad$ & $\dst \quad so(8)\oplus
so(3) \quad$ & $\dst \quad spin(7)\oplus
so(3) \quad$ &    \\
\Cline{1.2pt}{1-5}
\end{tabular}
\medskip\centerline{Table 1. Tamaru's classification \cite{Ta}}
\end{table}

\begin{theorem}
All pairs $(G,H)$ associated to the triples of Lie algebras
$(\mathfrak g,\mathfrak k,\mathfrak h)$ listed in Table 1 are
integrable. The geodesic flows on g.o. spaces $(G/H,ds^2_\lambda)$
are Liouville integrable by means of analytic integrals,
polynomial in momenta.
\end{theorem}

\noindent{\it Proof.} Integrability of the pairs $(G,H)$ in the
cases 1, 2, 9, 5, 8 is described, respectively, in Examples
\ref{flag1}, \ref{sfera1}, \ref{sfera2} and \ref{projektivni},
while in the cases 10, 11, 12 it follows from Theorem \ref{pomoc},
item (ii) and from the fact that symmetric pairs are integrable.
Further, we can simply use the known classifications of complexity
0 and 1 homogeneous spaces given in \cite{Mik, Ya, Pa, MS, ArCh}
to conclude that the others homogeneous spaces in Tamaru's
classification have the complexity equal 0 or 1 as well. Thus, all
pairs $(G,H)$ associated to the triples of Lie algebras
$(\mathfrak g,\mathfrak k,\mathfrak h)$ listed in Table 1 are
integrable. \hfill$\Box$

\begin{rem}{\rm
Up to the covering, which preserves the integrability of geodesic
flows, the g.o. spaces of compact simple Lie groups obtained in
\cite{KoVa, AlAr, Ta0} are also included in Table 1 (for example,
the weakly symmetric spaces $Spin(8)/G_2$ and $Spin(9)/Spin(7)$
\cite{Zi2, Ta0} correspond to the cases 3 and 4). On the other
hand, since $(SO(5)\times SO(2))/U(2)$ is a weakly symmetric space
\cite{Ya}, the g.o. spaces $((SO(5)\times SO(2))/U(2),ds^2_{p,q})$
constructed in \cite{DKN} have Liouville integrable geodesic flows
as well. }\end{rem}

\paragraph{4.2.}
Under the assumption of Theorem \ref{pomoc}, item (ii), if
$\dim\mathfrak l=1$ and the complexity of $G/K$ is equal to $0$,
we get that the complexity of the homogeneous space $G/H$ is
either equal to 0 (a linear function on $\mathfrak l$ is a central
function of $\R[\mathfrak v]^H$) or 1 (a linear function on
$\mathfrak l$ is not a central function), see also Remark 6 in
\cite{MS}.

For example, the flag manifolds \eqref{flag}, as weakly symmetric
spaces have the complexity $c=0$. Starting from $SO(2n+1)/U(n)$
and $Sp(n)/U(1)\times Sp(n-1)$ we obtain the complexity one
homogeneous spaces $SO(2n+1)/SU(n)$ and $Sp(n+1)/Sp(n)$. The
spaces $Sp(n)/U(1)\times Sp(n-1)$ and $Sp(n+1)/Sp(n)$ will be
treated in Examples \ref{projektivni} and \ref{sfera1} below. In
the next example we derive the complexity of homogeneous spaces
$SO(2n+1)/U(n)$ and $SO(2n+1)/SU(n)$.

\begin{exm}\label{flag1}{\rm
 Consider a triple of
Lie groups $(SO(2n+1),U(n),SU(n))$. The inclusion $u(n)\subset
so(2n+1)$ is the standard one, given by the mapping
$$
P+iQ \in u(n) \longmapsto \begin{pmatrix}
P & Q & 0\\
 -Q & P & 0\\
 0 & 0 & 0
\end{pmatrix}\in so(2n+1).
$$

Let us take $x_\mathfrak m=x_1+x_2$, where
\begin{eqnarray*}
&&x_1=\alpha_1 (E_1\wedge E_2-E_{2m+1}\wedge
E_{2m+2})+\dots+\alpha_m (E_{2m-1}\wedge E_{2m}-E_{4m-1}\wedge
E_{4m}), \\
&& x_2=\beta_1 E_1\wedge E_{4m+1}+\dots+\beta_{4m} E_1\wedge
E_{4m+1}, \qquad m=[n/2].
\end{eqnarray*}

Straightforward computations yield $u(n)_{x_1} \cong su(2)^m\oplus
u(1)^{n-2m}$, $u(n)_{x_2}\cong u(n-1)$ and
$$
u(n)_{x_\mathfrak m}=u(n)_{x_1}\cap u(n)_{x_2}=\{0\},
$$
for a generic values of the parameters $(\alpha_1,\dots,\alpha_m)$
and $(\beta_1,\dots,\beta_{4m})$. Thus, $\dim \mathfrak
u(n)_{x_\mathfrak m}=0$, for a generic $x\in\mathfrak m$.

Further, $\dim\mathfrak m-\dim u(n)=\dim so(2n+1)-2\dim
u(n)=n=\rank so(2n+1)$. Therefore, $\ddim\R[\mathfrak
m]^{U(n)}=\dind\R[\mathfrak m]^{U(n)}=n$ (see \eqref{ddim} and
\eqref{dind}) and $c(SO(2n+1),U(n))=0$. A complete set of
independent commuting $\Ad_{U(n)}$-invariant polynomials is simply
\begin{equation*}\label{flag2}
p_1=\tr(x^2),\,p_2=\tr(x^4),\,\dots,\,p_n=\tr(x^{2n}).
\end{equation*}
Next, according to item (ii) of Theorem \ref{pomoc}, a complete
commutative set on $\mathfrak v$ has $\dind\R[\mathfrak
m]^{U(n)}+1= \rank so(2n+1)+1$ independent functions. Thus,
$c(SO(2n+1),SU(n))=1$.
 }\end{exm}

\begin{exm}\label{sfera1}{\rm
{ $(G,K,H)=(Sp(n+1),Sp(1)\times Sp(n),Sp(n))$, $G/H=S^{4n+3}$.} We
have $\mathfrak v=\mathfrak l \oplus \mathfrak m \cong sp(1)\oplus
\mathbb H^n$. The $\Ad_{Sp(n)}$-action on the first factor is
trivial, while on  $\mathbb H^n=\R^{4n}$ is the standard one: the
orbit of $e\in\mathbb H^n$ is a sphere $S^{4n-1} \cong
Sp(n)/Sp(n-1)$ (see \cite{Zi}). Therefore, the algebra
$\R[\mathfrak v]^{Sp(n)}$ is generated by the linear functions on
$sp(1)$ and a quadratic invariant on $\mathbb H^n$, while the
center of $\R[\mathfrak v]^{Sp(n)}$ is generated by quadratic
invariants on $sp(1)$ and $\mathbb H^n$ (for example $\delta(x)$
and $h_0(x)$). Thus
$$\ddim\R[\mathfrak v]^{Sp(n)}=4, \qquad
\dind\R[\mathfrak v]^{Sp(n)}=2.
$$
For a complete commutative set $\mathcal B$ we can take an
arbitrary linear function $f$ on $sp(1)$, together with
$\delta(x)$ and $h_0(x)$. In particular we get the complexity
$c(Sp(n+1),Sp(n))=1$.

The sphere $(S^{4n+3}=Sp(n+1)/Sp(n),ds^2_\lambda)$ has
$Sp(n+1)\times Sp(1)$ as a full group of isometries (see the proof
of Theorem \ref{pomoc}, item (i)) and
$$
S^{4n+3}=Sp(n+1)\times Sp(1)/Sp(n)\times Sp(1)
$$ is a weakly
symmetric space, a distance sphere in a quaternionic projective
space $\mathbb{HP}^{n+1}=Sp(n+2)/Sp(n+1)\times Sp(1)$ (see
\cite{Zi, BeVa}). From the analysis above, we see that  the
commutative algebra of $Sp(n+1)\times Sp(n)$-invariant polynomials
on $T^*S^{4n+3}$ is generated by $\delta(x)$ and $h_0(x)$.
}\end{exm}

\begin{exm}\label{sfera2}{\rm
{\ $(G,K,H)=(SU(n+1),U(n),SU(n))$, $G/H=S^{2n+1}$.} Now,
$\mathfrak v=\mathfrak l \oplus \mathfrak m \cong \R \oplus
\mathbb C^n$. The $\Ad_{SU(n)}$-action on $\mathbb C^n=\R^{2n}$ is
the usual action and it is trivial on $\R$ (see \cite{Zi}).
Therefore, $ \ddim\R[\mathfrak v]^{SU(n)}=\dind\R[\mathfrak
v]^{SU(n)}=2 $, $\R[\mathfrak v]^{SU(n)}$ is generated by $h_0(x)$
and by a linear function on $\mathfrak l\cong\R$,  and
$c(SU(n+1),SU(n))=0$. Note that $
(S^{2n+1}=SU(n+1)/SU(n),ds^2_\lambda)$ can be seen as a distance
sphere in a complex projective space
$\mathbb{CP}^{n+1}=SU(n+2)/S(U(n+1)\times U(1))$ (see \cite{Zi,
BeVa}).}\end{exm}

\begin{exm}\label{projektivni}{\rm {$(G,K,H)=(Sp(n+1),Sp(1)\times
Sp(n),U(1)\times Sp(n))$, $G/H=\mathbb{CP}^{2n+1}$.} Here we have
a natural inclusion $U(1)\subset Sp(1)$, $sp(1)=u(1)\oplus \mathbb
C$ and $\mathfrak l\cong\mathbb C$, $\mathfrak m\cong \mathbb
H^n$. $U(1)$ acts by a rotation on $\mathbb C$ and trivially on
$\mathbb H^n$ and $Sp(n)$ acts trivially on $\mathbb C$ and by its
standard representation on $\mathbb H^n$. Thus $ \ddim\R[\mathfrak
v]^{Sp(n+1)}=\dind \R[\mathfrak v]^{Sp(n+1)}=2 $, $\R[\mathfrak
v]^{Sp(n+1)}$ is generated by $h_0(x)$ and $\delta(x)$, and
$c(Sp(n+1),U(1)\times Sp(n))=0$.}
\end{exm}

\subsection*{Acknowledgments} I am greatly thankful to
Dmitri Alekseevsky for pointing out the problem of geodesic flows
on g.o. spaces and to the referee for useful remarks which helped
me to improve the exposition and to extend some results of the
paper. This research was supported by the Serbian Ministry of
Science, Project 174020 Geometry and Topology of Manifolds,
Classical Mechanics and Integrable Dynamical Systems.


\begin{thebibliography}{84}
 \small

\bibitem{AlAr}
Alekseevsky, D. and Arvanitoyeorgos, A.: Riemannian flag manifolds
with homogeneous geodesics.  Trans. Amer. Math. Soc. {\bf 359}
(2007), no. 8, 3769–-3789

\bibitem{AlNi}
Alekseevsky, D. V.  and  Nikonorov, Y.G.: Compact Riemannian
manifolds with homogeneous geodesics.  SIGMA  {\bf 5}  (2009),
Paper 093, 16 pp, arXiv:0904.3592 [math.DG]

\bibitem{ADN}
{Arvanitoyeorgos, A., Dzhepko, V. V. and Nikonorov, Yu. G.}:
Invariant Einstein metrics on some homogeneous spaces of classical
Lie groups, Canadian Journal of Mathematics  {\bf 61}  (2009), no.
6, 1201--1213,  arXiv: math/0612504 [math.DG].


\bibitem{ArCh} Arzhantsev, I. V. and Chuvashova, O. V.:
Classification of affine homogeneous spaces of complexity one.
(Russian)  Mat. Sb.  {\bf 195}  (2004), no. 6, 3--20;  translation
in Sb. Math.  195  (2004),  no. 5-6, 765–-782

\bibitem{BeVa} Berndt, J. and Vanhecke, L.: Geometry of weakly symmetric spaces.
J. Math. Soc. Japan  {\bf 48}  (1996),  no. 4, 745–-760.

\bibitem{BKV}
Berndt, J., Kowalski, O. and Vanhecke, L.: Geodesics in weakly
symmetric spaces.  Ann. Global Anal. Geom.  {\bf 15}  (1997), no.
2, 153-–156.

\bibitem{Be} Besse, A.: {Einstein Manifolds}, Springer, A Series of Modern Surveys in Mathematics,
(1987)

\bibitem{BJ1} Bolsinov, A. V.  and Jovanovi\' c, B.:
{Integrable geodesic flows on homogeneous spa\-ces}, Matem.
Sbornik \textbf{192} (2001), no. 7, 21--40 (Russian); English
translation: Sb. Mat. \textbf{192} no. 7--8, 951--969, (2001)


\bibitem{BJ2} Bolsinov, A. V. and Jovanovi\' c, B.:
{Non-commutative integrability, moment map and geodesic flows}.
Annals of Global Analysis and Geometry {\bf 23}, no. 4, 305--322
(2003), arXiv: math-ph/0109031

\bibitem{BJ3} Bolsinov, A. V. and Jovanovi\' c, B.:
Complete involutive algebras of functions on
cotangent bundles of homogeneous spaces.
Mathematische Zeitschrift {\bf 246} no.  1-2, 213--236 (2004)

\bibitem{Br} Brailov, A. V.:
{Construction of complete integrable geodesic flows on compact
symmetric spaces}. Izv. Acad. Nauk SSSR, Ser. matem. {\bf 50},
no.2, 661-674  (1986) (Russian); English translation: {Math.
USSR-Izv.} {\bf 50}, no.4, 19-31 (1986)

\bibitem{DGJ}
Dragovi\'c, V., Gaji\'c, B.  and Jovanovi\' c, B.: {Singular
Manakov Flows and Geodesic Flows of Homogeneous Spaces of
$SO(n)$}, Transfomation Groups {\bf 14}, no. 3, 513--530 (2009),
arXiv: 0901.2444

\bibitem{DKN}
Du\v sek, Z., Kowalski, O. and Nik\v cevi\' c, S.: New examples of
Riemannian g.o. manifolds in dimension 7.  Differential Geom.
Appl.  {\bf 21}  (2004),  no. 1, 65-–78.


\bibitem{FJ} Fedorov, Yu. N. and Jovanovi\' c, B.:
Geodesic flows and Neumann systems on Stiefel varieties: geometry
and integrability, Mathematische Zeitschrift {\bf 270}(2012),
659–-698,   arXiv:1011.1835 [nlin.SI]


\bibitem{Go} Gordon, C. S.:
Homogeneous Riemannian manifolds whose geodesics are orbits.
Topics in geometry,  155–-174, Progr. Nonlinear Differential
Equations Appl., 20, Birkhäuser Boston, Boston, MA, 1996.

\bibitem{Je} Jensen, G.: Einstein metrics on principal
fiber bundles, J. Diff. Geom. {\bf 8} (1973), 599-614.

\bibitem{Jo0} Jovanovi\' c, B.: Symmetries and integrability.  Publ. Inst.
Math. (Beograd) (N.S.)  {\bf 84(98)}  (2008), 1–-36.
arXiv:0812.4398 [math.SG]

\bibitem{Jo}Jovanovi\' c, B.:
Integrability of Invariant Geodesic Flows on n-Symmetric Spaces,
Annals of Global Analysis and Geometry,  {\bf 38} (2010),
305--316, arXiv:1006.3693 [math.DG]

\bibitem{Ka} Kaplan, A.:
On the geometry of groups of Heisenberg type.  Bull. London Math.
Soc.  {\bf 15}  (1983), no. 1, 35–-42.

\bibitem{KoVa}
Kowalski, O. and Vanhecke, L.: Riemannian manifolds with
homogeneous geodesics. Boll. Un. Mat. Ital. B (7) 5 (1991), no. 1,
189–-246.

\bibitem{KoNi}
Kowalski, O. and Nik\v cevi\' c, S.: On geodesic graphs of
Riemannian g.o. spaces.  Arch. Math. (Basel)  {\bf 73}  (1999),
no. 3, 223–-234; Appendix:  Arch. Math. (Basel)  {\bf 79}  (2002),
no. 2, 158–-160.

\bibitem{La} Lacomba, E. A.:
Mechanical systems with symmetry on homogeneous spaces. Trans.
Amer. Math. Soc. {\bf 185} (1973), 477-–491.

\bibitem{LMV} Laurent-Gengoux, C., Miranda, E.  and Vanhaecke, P.:
Action-angle Coordinates for Integrable Systems on Poisson
Manifolds, International Mathematics Research Notices,
doi:10.1093/imrn/rnq130, arXiv: 0805.1679 [math.SG].

\bibitem{MF2}  Mishchenko, A. S. and Fomenko, A. T.:
{Generalized Liouville method of integration of Hamiltonian
systems}. Funkts. Anal. Prilozh. {\bf 12}, No.2, 46-56  (1978)
(Russian); English translation: Funct. Anal. Appl. {\bf 12},
113--121  (1978)

\bibitem{Mik} Mikityuk, I. V.:
{On the integrability of invariant Hamiltonian systems with
homogeneous configuration spaces.}  {Matem. Sbornik}, {\bf
129(171)}, no.4, 514-534 (1986) (Russian);   English translation:
{Math. USSR Sbornik}, {\bf 57}, no.2, 527-547 (1987)

\bibitem{MP} Mykytyuk, I. V.  and Panasyuk, A.:
{Bi-Poisson structures and integrability of geodesic flows on
homogeneous spaces}, Transformation Groups \textbf{9}(3),
289--308, (2004)

\bibitem{MS} Mykytyuk, I. V. and Stepin, A. M.:
{Classification of almost spherical pairs of compact simple Lie
groups}. In: {\it Poisson geometry}, Banach Center Publ., 51,
Polish Acad. Sci., Warsaw, 2000, pp. 231-241.


\bibitem{Mik2}
Mykytyuk, I. V.: Actions of Borel subgroups on homogeneous spaces
of reductive complex Lie groups and integrability. Composito Math.
{\bf 127}, 55-67, (2001).


\bibitem{N} Nekhoroshev, N. N.: Action-angle variables
and their generalization. Tr. Mosk. Mat. O.-va. {\bf 26},
181--198, (1972)  (Russian); English translation: {Trans. Mosc.
Math. Soc.} {\bf 26},  180--198 (1972)


\bibitem{Pa}
Panyushev, D. I.: Complexity of quasiaffine homogeneous varieties,
t--decompositions, and affine homogeneous spaces of complexity 1,
Adv. in Soviet Math. {\bf 8}, 151-166 (1992).

\bibitem{Sa} Sadetov, S. T.: A proof of the Mishchenko-Fomenko conjecture
(1981).  Dokl. Akad. Nauk {\bf 397}, no. 6, 751--754 (2004)
(Russian)

\bibitem{Se} Selberg, A.: Harmonic analysis and discontinuous groups in
weakly symmetric Riemannian spaces with applications to Dirichlet
series. J. Indian Math. Soc. (N.S.) {\bf 20} (1956), 47–-87.

\bibitem{Ta0}
Tamaru, H.: Riemannian geodesic orbit metrics on fiber bundles.
Algebras Groups Geom.  {\bf 15}  (1998),  no. 1, 55-–67.

\bibitem{Ta} Tamaru, H.:
Riemannian G.O. spaces fibered over irreducible symmetric spaces.
Osaka J. Math. {\bf 36} (1999), no. 4, 835-–851.

\bibitem{Th} Thimm A.:
{Integrable geodesic flows on homogeneous spaces}, Ergod. Th. \&
Dynam. Sys.,{\bf 1}, 495--517  (1981)

\bibitem{To} T\'oth, G. Z.:
On Lagrangian and Hamiltonian systems with homogeneous
trajectories, J. Phys. A: Math. Theor. {\bf 43} (2010) 385206,
arXiv: 1003.1495v3 [math-ph]


\bibitem{Vi} Vinberg, E. B.:
Commutative homogeneous spaces and co-isotropic symplectic
actions. Uspekhi Mat. Nauk {\bf 56}, No. 1, 3-62 (2001) (Russian);
English translation: Russian Math. Surveys {\bf 56}, No. 1, 1-60
(2001)

\bibitem{Ya}
Yakimova, O. S. Weakly symmetric Riemannian manifolds with a
reductive isometry group. (Russian)  Mat. Sb.  {\bf 195}  (2004),
no. 4, 143--160;  translation in  Sb. Math.  {\bf 195}  (2004),
no. 3--4, 599-–614

\bibitem{Zi} Ziller, W.: Homogeneous Einstein metrics on Spheres
and Projective Spaces, Math. Ann. {\bf 259} (1982) 351--358.

\bibitem{Zi2}
Ziller, W.: Weakly symmetric spaces.  Topics in geometry,
355-–368, Progr. Nonlinear Differential Equations Appl., {\bf 20},
Birkhäuser Boston, Boston, MA, 1996.

\end{thebibliography}
\end{document}